\documentclass{article}
\usepackage[utf8]{inputenc}
\usepackage[pass]{geometry}

\usepackage[utf8]{inputenc}
\usepackage[T1]{fontenc}
\usepackage{lmodern}

\usepackage{listings}

\usepackage{amsmath}
\usepackage{mathtools}
\usepackage{amssymb}
\usepackage{amsthm}

\usepackage{array}

\usepackage[plainpages=false,pdfpagelabels,colorlinks=true,citecolor=blue,hypertexnames=false]{hyperref}

\newtheorem{thm}{Theorem}

\newtheorem{lem}{Lemma}

\theoremstyle{definition}

\newtheoremstyle{note}
{3pt}
{3pt}
{}
{}
{\bfseries}
{.}
{.5em}
{}
\theoremstyle{note}

\newtheorem{ex}{Example}

\newmuskip\pFqmuskip

\newcommand*\pFq[6][8]{%
  \begingroup 
  \pFqmuskip=#1mu\relax
  \mathcode`\,=\string"8000
  \begingroup\lccode`\~=`\,
  \lowercase{\endgroup\let~}\pFqcomma
  {}_{#2}F_{#3}{\left[\genfrac..{0pt}{}{#4}{#5};#6\right]}%
  \endgroup
} 
\newcommand{\pFqcomma}{\mskip\pFqmuskip}

\newcommand{\ps}[1]{[\![#1]\!]}
\newcommand{\FF}[0]{\mathbb{F}}
\newcommand{\NN}[0]{\mathbb{N}}
\newcommand{\ZZ}[0]{\mathbb{Z}}
\newcommand{\QQ}[0]{\mathbb{Q}}

\newcommand{\gfun}{\textrm{\texttt{gfun }}}
\newcommand{\Iso}{\textrm{\sf Iso}}

\title{The art of algorithmic guessing in \texttt{gfun}}
\author{Sergey Yurkevich}

\begin{document}

\maketitle

\begin{abstract}
    The technique of guessing can be very fruitful when dealing with sequences which arise in practice. This holds true especially when guessing is performed algorithmically and efficiently. One highly useful tool for this purpose is the package named \texttt{gfun} in the software Maple. In this text we explore and explain some of \texttt{gfun}'s possibilities and illustrate them on two examples from recent mathematical research by the author and his collaborators. 
\end{abstract}

\section{Introduction}\label{sec:intro}

George Pólya described the mathematical scientific method by the short but precise premise ``First guess, then prove''~\cite[p.~27]{Polya78}. This mantra lies in the heart of experimental mathematics, which, while having many facets and branches, can be roughly described as a 3-step process: compute a high-order approximation of a problem, guess/conjecture a general pattern, prove the conjecture. Many years passed since Pólya formulated his advice and modern mathematicians can now profit of better hardware and efficient algorithms designed for not only proving statements but also guessing them.

Nowadays, experimental mathematics is almost impossible to imagine without the computational power provided to us by recent technology and theoretical algorithmic breakthroughs. Many principles in ``guessing'' and ``proving'' were adapted and can be run completely automatized taking just fractions of seconds in time, while finding and justifying highly non-trivial theorems. However, experience also shows that still many techniques which are self-evident to some scientists are inaccessible or unknown to others, even though these methods could often extensively aid their research. 

One class of objects which is ideal for the algorithmic ``guess and prove'' strategy turns out to be the class of holonomic sequences (to be defined below in~\S\ref{sec:precdfin}). This type of sequences is not only very common throughout mathematics and other sciences, but also provides a good algorithmic data structure and is very well studied from the theoretical point of view. Moreover, in Maple users can profit of an excellent package called \gfun designed for manipulating and guessing holonomic sequences efficiently. Similar packages exist for other computer algebra systems like, for example, in Mathematica (e.g. \href{https://www.mat.univie.ac.at/~kratt/rate/rate.html}{RATE}, see \cite{Krattenthaler99}, or \href{https://www3.risc.jku.at/research/combinat/software/ergosum/RISC/Guess.html}{Guess} by Kauers) or in SageMath (e.g. \href{http://fricas.sourceforge.net/}{Guess} in FriCAS by Rubey). 

The goal of this text is twofold. Our first aim is to give a short introduction to the practical use of Maple's package \texttt{gfun}, its various functions and the resulting possibilities. Secondly, we wish to explain the ``art of algorithmic guessing'' by means of two recent and very different results in mathematical research which were deduced with \texttt{gfun}'s help. The ultimate target is to convince the reader that algorithmic guessing is very powerful and at the same time easy to do efficiently in Maple. 

The structure of the paper is as follows: in the Introduction (\S\ref{sec:intro}) we explain the preliminaries and at the same time provide motivation for the objects of interest. We will describe very briefly the theoretical parts of the ``guess and prove'' technique in this context. In \S\ref{sec:gfun} we concentrate on illustrating the practical use of Maple's package \gfun by describing the most useful functions and illustrating them on some ``toy examples''. Finally, Section~\ref{sec:examples} is devoted to two recent applications from mathematical research by the author and his collaborators, which nicely demonstrate the power of guessing as well as the potential of \texttt{gfun}.

\subsection{P-recursive sequences and D-finite functions}\label{sec:precdfin}
In this text we will mostly work with the class of P-recursive sequences. This notion became already quite classical in computer science and mathematics, but we nevertheless recall that a sequence $(u_n)_{n\geq0}$ is called \emph{P-recursive} (or \emph{holonomic}) if it satisfies a linear recurrence with polynomial coefficients:
\begin{equation}\label{eq:prec}
    c_r(j) u_{j+r} + \dots + c_0(j) u_{j} = 0 , \quad j \geq 0, \quad c_r(x) \neq 0.
\end{equation}
P-recursive sequences are ubiquitous in mathematics and often also appear in other disciplines like physics or biology. Some prominent examples are: 
\begin{itemize}
    \item The \href{https://oeis.org/A000045}{Fibonacci sequence} $(F_n)_{n\geq0} = (0,1,1,2,3,5,\dots)$ which satisfies 
    \[
    F_{n+2} - F_{n+1} - F_{n} = 0.
    \]
    \item The sequence of \href{https://oeis.org/A000108}{Catalan numbers} $(C_n)_{n\geq0} = (1, 1, 2, 5, 14, 42,\dots)$:
    \[
    (n+2)C_{n+1} - 2(2n+1)C_n = 0.
    \]
    \item The \href{https://oeis.org/A000142}{factorial sequence} $(n!)_{n\geq0} = (1,1,2,6,24,\dots)$:
    \[
    (n+1)! - (n+1)n! = 0.
    \]
    \item The \href{https://oeis.org/A005259}{Apéry numbers} $(A_n)_{n \geq 0} = (1, 5, 73, 1445, \dots)$ which played a crucial role in the proof of the irrationality of $\zeta(3)$~\cite{Apery79, vdPoorten79}:
    \[
    (n+2)^3A_{n+2} - (2n + 3)(17n^2 + 51n + 39)A_{n+1} + (n+1)^3 A_{n} = 0.
    \]
    \item The Yang-Zagier numbers $(a_n)_{n\geq0} = (1, -48300, 7981725900,\dots)$ which appeared first in~\cite[p.~768]{Zagier18} and are investigated in \S\ref{sec:zag}.
    \item The numbers $(d_n)_{n\geq0} = (72, 1932, 31248,790101/2,\dots)$ which play an important role in the study of Canham's model in biology, see \S\ref{sec:iso}.
\end{itemize}
Like in the examples above, for the sake of simplicity and explicitness we will assume that $u_n \in \QQ$ for a holonomic sequence $(u_n)_{n\geq0}$. Of course, the definitions and almost all properties addressed below translate to arbitrary fields (usually, but not necessarily, of characteristic 0).

One reason for the importance of P-recursive sequences is the equivalent characterization on the level of generating functions. We recall that a formal power series $f(x) \in \QQ\ps{x}$ is called \emph{D-finite} (or \emph{holonomic}) if it satisfies a linear differential equation with polynomial coefficients:
\begin{equation}\label{eq:dfin}
    p_s(x)f^{(s)}(x) + \cdots + p_1(x) f'(x)+ p_0(x)f(x) = 0, \quad p_s(x) \neq 0.
\end{equation}
The proof of the following classical theorem connecting P-recursive sequences and D-finite functions can be found in Stanley's seminal article~\cite{Stanley80} (Theorem~1.5) where he mentions that Jungen~\cite{Jungen31} was already using it almost half a century before.
\begin{thm}\label{thm:recdfin}
    A sequence $(u_n)_{n\geq0}$ is P-recursive if and only if the generating function $\sum_{n\geq0} u_n x^n$ is D-finite.
\end{thm}
\noindent
General interest for the notion of holonomicity mostly comes from two main aspects:

\subsubsection{Algorithmic and mathematical theory}

From the algorithmic point of view, holonomic objects are useful because they only require finitely many data to be stored uniquely: on the level of P-recursive sequences it is clear that the polynomials $c_0(x), \dots, c_r(x)$ together with the initial terms $u_0,\dots,u_N$, where $N$ is the maximum between $r$ and the largest integral root of $c_r(x)$, are enough to encode a given sequence. The general fact that linear differential equations form an excellent data structure from the computational point of view is propagated by Salvy, for example in~\cite{Salvy19}.

P-recursive sequences/D-finite functions enjoy nice and effective closure properties. If $(u_n)_{n\geq0}$, $(v_n)_{n\geq0}$, $f(x)$ and $g(x)$ are holonomic, then
\begin{enumerate}
    \item $(u_n+v_n)_{n\geq0}$ and $f(x)+g(x)$ are holonomic.
    \item $(u_n \cdot v_n)_{n\geq0}$ and $f(x)\cdot g(x)$ are holonomic.
    \item $f'(x)$ and $\int f(x) \mathrm{d}x$ are holonomic.
\end{enumerate}
For the proofs of all these properties we refer to \cite[Thm. 2.3 \& Thm. 2.10]{Stanley80}. It essentially only uses linear algebra and the property that $f$ is D-finite if and only if $f,f',f'',\dots$ span a finite-dimensional vector space over $\QQ(x)$. Moreover, a theorem sometimes attributed to Abel states that:
\begin{enumerate}
    \setcounter{enumi}{3}
    \item If $a(x)$ is an algebraic function\footnote{Recall that $a(x) \in \QQ\ps{x}$ is called \emph{algebraic} if there exists a bivariate non-zero polynomial $P(x,y)$ in $\QQ[x,y]$ such that $P(x,a(x)) = 0$. A non-algebraic series is called \emph{transcendental}.} then $a(x)$ is holonomic.
\end{enumerate}
More generally, if $a(x)$ is algebraic and $f(x)$ D-finite then~\cite[Thm.~2.7]{Stanley80}
\begin{enumerate}
    \setcounter{enumi}{4}
    \item $f(a(x))$ is holonomic as well.
\end{enumerate}
All these closure properties are effective in the sense that there exist algorithms for performing them on the level of differential equations and recursions. As we will see, they are implemented efficiently in Maple's package \texttt{gfun}. 

We warn the reader that the quotient and composition of holonomic functions is not necessarily holonomic. For example the functions $\tan(x)$ and $e^{e^x-1}$ are not D-finite. Similarly, $(1/F_n)_{n\geq1}$, where $F_n$ is the $n$-th Fibonacci number is not P-recursive.

From the theoretical point of view, linear differential equations have been studied by many outstanding mathematicians in the past two centuries. Since giving a short and at the same time complete overview is impossible, we refer the interested reader to some of the best introductory books on this topic for a huge universe of amazing theorems and discoveries~\cite{Poole60, Gray00, PuSi03}. 

A linear differential equation like (\ref{eq:dfin}) is usually studied from the viewpoint of the attached \emph{differential operator} 
\[
    p_s(x) \partial^s + \cdots + p_1(x)\partial + p_0(x) \in \QQ[x]\langle \partial \rangle,
\]
where the symbol $\partial$ stands for $\frac{\mathrm{d}}{\mathrm{d}x}$ and we have the commutation rule $\partial x = x \partial +1$. The (non-commutative) polynomial ring $\QQ[x]\langle \partial \rangle$ is called the \emph{Weyl algebra}. The order of an operator $L \in \QQ[x]\langle \partial \rangle$ is defined as the largest degree of $\partial$. A classical result, usually attributed to Ore~\cite{Ore32}, states that $\QQ(x)\langle \partial \rangle$ is a Euclidean domain. Therefore it makes sense to speak about the so-called GCRD (greatest common right divisor) and LCLM (least common left multiple) of two operators $A,B \in \QQ(x)\langle \partial \rangle$. Naturally, LCLM$(A,B)$ is defined as the least-order monic operator $L \in \QQ(x)\langle \partial \rangle$ such that $L=Q_1 A=Q_2B$ for some non-zero $Q_1,Q_2 \in \QQ(x)\langle \partial \rangle$, and $G=\mathrm{GCRD}(A,B)$ is the monic operator of largest order such that $A = Q_1G$ and $B=Q_2G$ for some $Q_1,Q_2 \in \QQ(x)\langle \partial \rangle$. It can be proved that the solution spaces of $Ay=0$ and $By=0$ are vector spaces and that the LCLM produces an operator whose solution space is the sum of the spaces, while the solution space of GCRD$(A,B)$ is their intersection. 

Moreover, the LCLM and GCRD can be algorithmically computed using the ``skew version'' of the Euclidean algorithm; in practice, for example, with Maple's package \texttt{DEtools}. Then some elementary reasoning implies that the question of equality of two D-finite functions (and consequently two P-recursive sequences) is decidable and usually easy to answer in practice, see~\cite[\S4.3]{BoLaSa17}. We will elaborate on this in Section \ref{sec:examples} on practical examples.

Still on the algorithmic side, we can rely on many relatively recent and outstanding works by many still active scientists. For example, highly efficient algorithms for computing analytic continuation and evaluation of D-finite functions, based on ideas of the Chudnovsky brothers~\cite{ChCh90} were created by van der Hoeven~\cite{Hoeven99, Hoeven01, Hoeven07} and Mezzarobba~\cite{Mezzarobba10}. Notably, as we will see later, very useful in practice is the work by van Hoeij and collaborators~\cite{KuHo13,HoVi15,ImHo15} on explicit solutions of differential equations in terms of known special functions. The famous method of \emph{creative telescoping} propagated by Zeilberger~\cite{Zeilberger91}, and constantly improved in the last decades, allows for finding and proving linear recurrences/differential equations for sequences given as explicit sums or functions given as integrals; see for example~\cite{Chyzak14} for a great exposition of many achievements in this field and \cite{ChKa17} for open problems. We also mention very recent works by Bostan, Rivoal, Salvy~\cite{BoRiSa21} and by Barkatou, Cluzeau, Di Vizio, Weil~\cite{BaClDiWe20} which allow for practical proofs of transcendence and algebraicity of given D-functions.

Finally, as already mentioned, P-recursive sequences and D-finite functions happen to form a great class for efficient guessing algorithms. We will elaborate on this in~\S\ref{sec:guessnprove}.

\subsubsection{Importance in practice}
Interest in holonomic objects is also motivated by the fact that, as experience shows, they often appear in practice. Their practical importance can already be guessed from the huge amount of existing theory and algorithms. Closure properties show how easy it is to build arbitrarily complicated examples of D-finite functions and P-recursive sequences. 

As a further example, we recall the notion of the \emph{Gaussian hypergeometric function} ${}_2F_1$ defined by
\begin{equation}\label{eq:2F1}
\pFq{2}{1}{a,b}{c}{x} \coloneqq \sum_{n=0}^\infty \frac{(a)_n(b)_n}{(c)_n} \, \frac
{x^n} {n!},
\end{equation}
where $(u)_j \coloneqq u(u+1)\cdots(u+j-1)$ is the rising factorial and $a,b\in \QQ$, $c \in \QQ \setminus \ZZ_{<0}$ are parameters. It is not difficult to see that the coefficient sequence of such a function satisfies a first-order linear recurrence relation with polynomial coefficients, whereas the function itself satisfies a linear differential equation of order~2. These generating functions are first non-trivial examples from the holonomic viewpoint, however they already play an important role in enumerative combinatorics, as well as in the study of elliptic functions and special functions in general, modular forms, orthogonal polynomials, etc. We will see them appearing in completely different contexts in \S\ref{sec:zag} and \S\ref{sec:iso}.

There exist (informal) estimates that the proportion of sequences in the \href{https://oeis.org/The On-Line Encyclopedia of Integer Sequences}{The On-Line Encyclopedia of Integer Sequences}~\cite{OEIS} which are P-recursive is around $25\%$; moreover, roughly $60\%$ of functions are D-finite in the Handbook of Mathematical Functions~\cite{AbSt64}\footnote{These numbers were estimated by Salvy~\cite{Salvy05} in 2005 and the magnitude of the first proportion was confirmed by the author as of 2022.}. These numbers, however, do not say much without closer inspection towards applications. Therefore, in order to support the thesis that holonomic objects appear in practical applications by means of an example, we will briefly recapitulate on the famous success story in enumerative combinatorics about the classification of lattice walks with small steps in the quarter plane. 

A \emph{walk} of length $n$ in the plane is a sequence $p_0,\dots,p_n$ of elements in $\ZZ^2$ such that $p_{i+1}-p_i \in \mathfrak{S}$, where $\mathfrak{S}$ is a finite subset of vectors in $\ZZ^2$. If all elements of $\mathfrak{S}$ have Euclidean length of at most $\sqrt{2}$ (i.e., $\mathfrak{S}$ consists only of vectors directing to the neighbours in $\ZZ^2$) one speaks of a \emph{small-step walk}. The walk is said to be \emph{confined to the upper half-space} or \emph{confined to the quarter plane} if moreover $p_0 = (0,0)$ and each $p_i \in \ZZ\times\NN$ in the first case or each $p_i \in \NN^2$ in the second case. Given a set $\mathfrak{S}$, a natural question is: how many unrestricted walks, or walks confined to the upper half-space or confined to the quarter plane of length $n$ exist? In particular, the nature of the generating functions is intriguing. 

It is easy to see that if there are no restrictions on the confinement, the generating function of the number of walks in the plane is rational. The case of short-step walks restricted to the half-space is more interesting: a result due to Bousquet-Mélou and Petkovšek says that the length generating function of such a walk is necessarily algebraic, see \cite{BoPe00,BoPe03}. The case of walks confined in a quarter plane is more difficult: it turns out that after disregarding uninteresting situations (e.g. $\mathfrak{S} = \{(-1,-1)\}$ where the number of walks is 0 for $n>0$), and after taking into account symmetry, exactly 79 cases remain. A great classification effort by many researchers was undertaken in the last decades. In particular, it is proved that exactly 23 of the 79 models (roughly 29\%) have a D-finite generating function and of those exactly 6 are algebraic. We refer to Bostan's habilitation thesis for an excellent summary on this topic~\cite{Bostan17}, see also \cite{Bostan21} for a short but clear-cut exposition.

\subsection{``Guess and prove'' for holonomic objects}\label{sec:guessnprove}
The guessing strategy for D-finite functions and P-recursive seqeunces is wonderfully explained in \cite[\S2]{Bostan17} and efficiently implemented in Maple's \texttt{gfun}. We will briefly summarize here the main theoretical ideas before we explain the practical use in Maple in the next section. 

Since our main object of interest are sequences, we will stick to the P-recursive viewpoint; in the case of linear differential equations instead of recursions, everything works analogously. 

The setting is the following: assume we are given some terms $u_0,\dots,u_N$ of a sequence $(u_n)_{n\geq0}$. In general we have no reason to assume that this sequence is P-recursive and we have no bounds on the order or degree of a possible recurrence. Still, we would like to guess a linear recurrence relation for $(u_n)_{n\geq0}$ from the data we have; if we are ``lucky'', we can prove afterwards that this recursion is indeed correct. The idea is to look for some natural number $r$ and polynomials $c_0(n),\dots,c_r(n)$ of some degree, say $d$, such that
\begin{equation}\label{eq:guessprec}
    c_r(j) u_{j+r} + \dots + c_0(j) u_{j} = 0
\end{equation}
holds for $j=0,\dots,N-r$. Clearly, if $(u_n)_{n\geq0}$ is P-recursive,  then such an equation exists for \emph{all}, arbitrarily large, natural numbers $N$, and for some, but fixed, $r,d \in \NN$. In general it is easy to see that the task boils down to solving a system of linear equations, where the unknowns are the $(r+1)(d+1)$ coefficients of the polynomials, see Example~\ref{ex:rec} below. Then, if $(r+1)(d+1)>N - r$, a linear algebra argument implies that a non-zero solution to this problem always exists. On the other hand, if $(r+1)(d+1) +r \ll N$, the system is clearly highly over-determined, and there is \textit{a priori} no reason for this system to have a solution, except, of course, if the sequence $(u_n)_{n\geq0}$ is P-recursive. Then it is also quite likely that the found solution is the true linear recurrence for $(u_n)_{n\geq0}$.

\begin{ex}\label{ex:rec}
Assume we are given the numbers 
\[
(u_n)_{0\leq n\leq6} = (1, 4, 36, 400, 4900, 63504,853776)
\] 
and we wish to find a linear recurrence of order one and degree at most two. In other words, we look for a non-zero sextuple $(a,b,c,d,e,f) \in \QQ^6$ such that
\[
(fn^2 + en + d) u_{n+1} + (cn^2 +bn+ a) u_n = 0, \quad n=0,\dots,5.
\]
Writing down what this means shows that we need to solve
\begin{align*}
    \begin{pmatrix}
        u_0 & 0 & 0 & u_1 & 0 & 0\\
        u_1 & u_1 & u_1 & u_2 & u_2 & u_2\\
        u_2 & 2u_2 & 4u_2 & u_3 & 2u_3 & 4u_3\\
        u_3 & 3u_3 & 9u_3 & u_4 & 3u_4 & 9u_4\\
        u_4 & 4u_4 & 16u_4 & u_5 & 4u_5 & 16u_5\\
        u_5 & 5u_5 & 25u_5 & u_6 & 5u_6 & 25u_6
    \end{pmatrix}
    \cdot
    \begin{pmatrix}
        a \\ b \\ c \\ d \\ e \\ f
    \end{pmatrix}
    =0.
\end{align*}
Indeed, using Gaussian elimination, we can easily find that the kernel of this matrix is spanned by $(-4,-16,-16,1,2,1)^t$ and therefore a valid candidate for our recursion in this example would be 
\[
(n^2 + 2n + 1)u_{n + 1} - (16n^2 + 16n + 4)u_{n} = 0. 
\]
\end{ex}

\bigskip

In practice there are usually no bounds on the order or degree of the recurrence one wants to guess. Therefore, the algorithm will first try $r=1$ and then successively increase the order $r$, while keeping $d \approx N/r$ such that the linear system is slightly over-determined. After a valid candidate for the recurrence is found, it is usually checked on a few known, but previously unused, terms of the sequence. \\

A very similar approach can be carried out when guessing the differential operator for a D-finite series or an annihilating polynomial for an algebraic function: in all these cases the problem reduces to linear algebra. In a specific but quite common situation in practice, such a guessing procedure almost immediately gives a proof. Let $f(x)$ be given as the solution of a linear ODE with enough initial conditions, and assume that we want to prove that $f(x)$ is algebraic. Assume, moreover, that we have a guess for an annihilating polynomial $P(x,y)$. Then proving its correctness is not difficult: using the effective version of the fact that algebraic functions are D-finite, we can convert $P(x,y)$ into a differential equation satisfied by all its roots and call $\Tilde{f}(x)$ its solution which coincides with $f(x)$ up to enough precision, such that $\Tilde{f}(x)$ is unique. Then the problem reduces to deciding equality of two D-finite functions which is algorithmically easy, see \S\ref{sec:zag} below for an example, or~\cite[\S4.3]{BoLaSa17} for the complete theoretical procedure. Exactly this strategy (in a slightly more general setting) was employed by Bostan and Kauers in their proof of the algebraicity of the generating function of Gessel walks~\cite{BoKa10}. \\

An important and highly non-trivial question is: How to guess \emph{efficiently}? In practice, as well as in the recent version of \texttt{gfun}, two clever ideas are applied in combination:
\begin{enumerate}
    \item Instead of solving the system over $\QQ$ it is faster to solve several systems over the finite fields $\FF_p$ for various prime numbers $p$ and then combine the solutions using the Chinese remainder theorem. In practice this method gives a huge speed-up.
    \item The linear systems one obtains from holonomic or algebraic guessing have structure that can exploited algorithmically. In fact, the problem can be viewed as a variant of Hermite-Padé approximation, for all types of which very efficient algorithms were invented by Beckermann and Labahn~\cite{BeLa97}. 
\end{enumerate}

\section{Working with \texttt{gfun}} \label{sec:gfun}
Before we explain and show the most useful functions of \texttt{gfun}, let us mention a bit of history about this package and credit the developers. 

The Maple package \gfun was written in 1992 by Salvy and Zimmermann. In the same year the authors submitted the article~\cite{SaZi94} which was intended to be an introduction to the package and the reference manual at the same time. Since then \gfun has been constantly improved by Salvy and got significantly better. Currently the newest version is 3.84 and can be downloaded from his web page\footnote{\href{http://perso.ens-lyon.fr/bruno.salvy/software/the-gfun-package/}{www.perso.ens-lyon.fr/bruno.salvy/software/the-gfun-package/}}. The package also comes pre-installed with Maple, however with the version 3.20 which is heavily outdated. One of the biggest improvements of the newest version compared to 3.20 is the efficiency boosts briefly explained on the theoretical level at the end of the previous section. Below we will indicate some more differences based on which we strongly recommend to download and use the newest version. Finally, we also mention the existence of the more recent Maple package \texttt{NumGfun} by Mezzarobba~\cite{Mezzarobba10} which deals with rigorous and efficient numerical computations with D-finite functions, since for some applications the two packages complement each other well. \\

Almost all functions of \gfun can roughly be divided into three categories:
\begin{itemize}
    \item Functions allowing interplay between the holonomic objects.
    \item Effective closure properties for P-recursive sequences/D-finite functions.
    \item Guessing algorithms for linear recurrences, ODEs and algebraic functions.
\end{itemize}
We will briefly summarize the most useful functions and showcase them on some ``toy examples''. We note that this section is quite similar to the original technical report on \gfun\cite{SaZi94}. However, as explained, in the last 30 years some things have changed and got significantly improved so we will highlight on them. \\

\noindent
\textbf{\texttt{listtoseries} and \texttt{seriestolist}:} These useful commands are trivial from the mathematical viewpoint: they transform a list of elements into a generating function (represented by a series) and vice versa. For example 
\begin{align*}
    &\texttt{> listtoseries([1, 1, 2, 4, 9],x);}\\
    &\texttt{$1 +x^2 + 2x^3 + 4x^4 + 9x^5 + O(x^6)$ }\\ 
    &\texttt{> seriestolist(1 +x\^{}2 + 2x\^{}3 + 4x\^{}4 + 9x\^{}5 + O(x\^{}6));}\\
    & [1,1,2,4,9]
\end{align*}

\noindent
\textbf{\texttt{rectodiffeq} and \texttt{diffeqtorec}} are the effective versions of Theorem~\ref{thm:recdfin}: they (efficiently and rigorously) convert a linear recurrence into a linear differential equation and vice versa:
\begin{align*}
    &\texttt{> rec := \{(n + 1)\^{}2 * u(n + 1) - 4*(2*n + 1)\^{}2 * u(n), u(0) = 1\}:}\\
    &\texttt{> deq := rectodiffeq(rec,u(n),y(x));}\\
    &\left\{-4y(x) + (-32x + 1)\frac{\mathrm{d}}{\mathrm{d}x}y(x) + (-16x^2 + x)\frac{\mathrm{d}^2}{\mathrm{d}x^2}y(x), \; y(0) = 1 \right\}\\
    &\texttt{> diffeqtorec(deq,y(x),u(n));}\\
    &\left\{(-16n^2 - 16n - 4)u(n) + (n^2 + 2n + 1)u(n + 1), u(0) = 1\right\}
\end{align*}

\noindent
\textbf{\texttt{`rec+rec`}, \texttt{`rec*rec`}, \texttt{`diffeq+diffeq`} and \texttt{`diffeq*diffeq`}:} These are implementations of the sum and product closure properties for P-recursive sequences and D-finite functions. For example, a linear recurrence relation for $(n! + 1/(n+1))_{n\geq0}$ can be found using
\begin{align*}
    &\texttt{> rec1 := \{u(n + 1) = u(n)*(n + 1), u(0) = 1\}:}\\
    &\texttt{> rec2 := \{{(n + 2)*u(n + 1) = (n + 1)*u(n), u(0) = 1}\}:}\\
    &\texttt{> `rec+rec`(rec1, rec2, u(n));}\\
    & (n + 1)(n + 3)u_{n+2} - (n + 2)(n^2 + 5n + 5)u_{n+1} + (n + 1)(n + 2)^2 u_n= 0.
\end{align*}
More general and in the same spirit are the commands \textbf{\texttt{poltorec}} and \textbf{\texttt{poltodiffeq}}.\\ 

\noindent
\textbf{\texttt{algeqtodiffeq}} is an effective implementation of Abel's theorem that every algebraic function is D-finite, see \S\ref{sec:precdfin}. It uses algorithmic ideas by Comtet~\cite{Comtet64} and the Chudnovsky brothers~\cite{ChCh86}, and is very efficient in practice. Clearly, together with \texttt{listtoalgeq} it is the central function for the algebraic guess-and-prove strategy, see \S\ref{sec:guessnprove}. \\

\noindent
\textbf{\texttt{holexprtodiffeq}} uses closure properties and known linear differential equations of special functions to output an ODE satisfied by the input function. For example:
\[
\texttt{> holexprtodiffeq(exp(sqrt(1 + x) - 1),y(x));}
\]
proves that $e^{\sqrt{1+x}-1}$ satisfies the differential equation
\[
-y(x) + 2 \frac{\mathrm{d}}{\mathrm{d}x}y(x) + (4x+4)\frac{\mathrm{d}^2}{\mathrm{d}x^2}y(x) = 0.
\]
We warn the user that the output is not necessarily the minimal-order differential equation. For example the output of 
\[
\texttt{> holexprtodiffeq(sin(x)cos(x),y(x));}
\]
is a third-order equation, whereas $\sin(x)\cos(x) = \sin(2x)/2$ clearly satisfies a linear second-order differential equation with polynomial coefficients as well.
\\

\noindent
\textbf{\texttt{rectohomrec} and \texttt{diffeqtohomdiffeq}:}
Holonomic objects are defined as solutions to linear and homogeneous equations with polynomial coefficients like (\ref{eq:prec}) and (\ref{eq:dfin}). It is not difficult to see that omitting the condition on homogeneity gives rise to the same class. These \gfun procedures allow to convert non-homogeneous equations to homogeneous ones.\\

\noindent
\textbf{\texttt{listtorec} and \texttt{seriestorec}} are the implementations of holonomic guessing for P-recursive sequences, which was briefly explained in \S\ref{sec:guessnprove}. As already elaborated, the efficiency of these algorithms has been significantly improved in the current version 3.84 compared to Maple's pre-installed \gfun 3.20. Moreover, thanks to the improved efficiency we can also profit of another significant upgrade: by default the old version guessed only linear recursions of order less or equal than 3 and whose coefficients have degree not larger than 4. While this is enough for the first four examples given in \S\ref{sec:precdfin}, many sequences in practice have larger defining equations. Admittedly, in the old version both of these figures could be changed by the user by specifying \texttt{gfun['maxordereqn']} and \texttt{gfun['maxdegcoeff']}. The version 3.84 dynamically adjusts these numbers ensuring maximal efficiency while maintaining the confidence in the guess.

With six terms we can guess the correct recurrence for the Catalan numbers. The input for \texttt{seriestorec} is a series, while for \texttt{listtorec} it should be a list of the elements.
\begin{align*}
    &\texttt{> ser := series((1 - sqrt(1 - 4*x))/(2*x),x,6);}\\
    & 1 + x + 2x^2 + 5x^3 + 14x^4 + 42x^5 + O(x^6)\\
    &\texttt{> seriestorec(ser, C(n));}
\end{align*}
finds $(-4n - 2)C_n + (n + 2)C_{n+1} = 0$. We need at least 15 terms of the \href{https://oeis.org/A125143}{Almkvist-Zudilin sequence} to find the correct recursion:
\begin{align*}
    &\texttt{> l := [1,-3,9,-3,-279,2997,-19431,65853,292329,-7202523,} \\
    & \qquad \texttt{69363009,-407637387,702049401,17222388453,-261933431751]:}\\
    &\texttt{> listtorec(l, u(n));}
\end{align*}
finds $(n + 2)^3u_{n+2} + (2n + 3)(7n^2 + 21n + 17)u_{n+1} + 81(n + 1)^3 u_n = 0.$

Note that if the guessing algorithm succeeded, the output of \texttt{listtorec} and \texttt{seriestorec} will actually be a list of two elements. The first is the guessed recurrence and the second is a flag variable: either \texttt{'ogf'} or \texttt{'egf'} depending on whether the algorithm could find a candidate for the ordinary sequence $(u_n)_{n\geq0}$ or for the exponential one $(u_n/n!)_{n\geq0}$. By indicating the \texttt{typelist} variable as a third input to the \gfun functions the user can force the algorithm to stick to one sequence type.  

We note that in the newest version of Maple 2021 great improvement on the package \texttt{LREtools} has been done by van Hoeij. One of the new available functions is \texttt{GuessRecurrence} which works essentially like \texttt{listtorec}.\\ 

\noindent
\textbf{\texttt{listtodiffeq} and \texttt{seriestodiffeq}} are commands analogous to the P-recursive guessing but in the D-finite case. Consequently, the above remarks on \gfun improvements also hold for these functions. The usage is the same, except that now the output is the guessed differential equation satisfied by the generating function of the input sequence. For example, the generating function of the squares of the central binomial coefficients satisfies:
\begin{align*}
    &\texttt{> l := [seq(binomial(2n,n)\^{}2,n=0..10)]:}\\
    &\texttt{> listtodiffeq(l, y(x));}\\
    &4y(x) + (32x - 1)\frac{\mathrm{d}}{\mathrm{d}x}y(x) + (16x^2 - x)\frac{\mathrm{d}^2}{\mathrm{d}x^2}y(x) = 0.
\end{align*}

\noindent
\textbf{\texttt{listtoalgeq} and \texttt{seriestoalgeq}:} These two functions are implementations of algebraic guessing. They guess an annihilating polynomial for the generating function of the input sequence. Similarly to the algorithms above, in \gfun 3.84 these functions are coded efficiently, i.e. they use modular arithmetic and fast algorithms for the resulting problems in linear algebra. For example, we find the minimal polynomial for the Motzkin numbers from the first 9 terms:
\begin{align*}
    &\texttt{> listtoalgeq([1, 1, 2, 4, 9, 21, 51, 127, 323], y(x)):}\\
    & 1 + (x - 1)y(x) + x^2y(x)^2 = 0.
\end{align*}

\noindent
\textbf{\texttt{rectoproc}} is a very useful \gfun function which is quite different from all the above ones. It translates a recurrence relation into a Maple procedure and allows to compute sequence elements conveniently and efficiently. If only the $N$-th term of a P-recursive sequence $(u_n)_{n\geq0}$  is needed, \texttt{rectoproc} uses binary splitting and fast integer multiplication which allows to find $u_N$ in quasi-optimal complexity. For example, let us compute the $N=10^5$-th Apéry number $A_N = \sum_{k=0}^N \binom{n+k}{k}^2 \binom{n}{k}^2$:
\begin{align*}
    &\texttt{> rec := \{(n+2)\^{}3*A(n+2) - (2*n + 3)*(17n\^{}2 + 51*n + 39)*A(n+1)} \\
    &\qquad \qquad \texttt{ + (n+1)\^{}3 A(n), A(0)=1, A(1)=5\}: }\\
    &\texttt{> pro := rectoproc(rec, A(n)):}\\
    &\texttt{> pro(10\^{}5):}
\end{align*}
On a regular laptop this finds $A_N$ in under 3 seconds. This is quite impressive as $A_N$ has more than 153 thousand digits. 

In practice one is often interested in the list of the first $N$ values of a sequence and not just in the $N$-th term. In this case \texttt{rectoproc(rec, u(n), list)} may be used.

\section{Two examples: recent applications}\label{sec:examples}
In this section we will focus on two examples from recent mathematical research which nicely demonstrate the usage of \gfun in a practical sense. We will briefly introduce the problems, then present and explain the algorithmic solutions. At the end of this section there is a short summary and conclusion containing the main takeaways of these examples.\\ 

Our first example comes from a recent survey article by Don Zagier~\cite{Zagier18} in which he investigates the structure of the generating function assigned to a particular integer sequence. The study of it and similar sequences is current work in progress by A. Bostan, J.-A. Weil and the author. 

The second example originates in a biological model describing the shape of biomembranes, the examination of which in a particular case culminates in the investigation of a specific function given by the quotient of two D-finite functions. The necessary study of it was first performed by Melczer and Mezzarobba~\cite{MeMe20} using rigorous asymptotics and then by Bostan and the author~\cite{BoYu21} from a completely different perspective.

\subsection{The Yang-Zagier numbers $(a_n)_{n\geq0}$}\label{sec:zag}
In \cite[p.~768]{Zagier18} Zagier introduces the integer sequence $(a_n)_{n\geq0}$, which he describes as a ``very mysterious example of numbers''. Following Zagier, we first define the P-recursive sequence $(c_n)_{n\geq 0}$ by
\begin{align*}
80352000 n(5 n-1)(5 n-2)(5 n-4) c_{n} & +\\ 
25\left(2592000 n^{4}-16588800 n^{3}+39118320 n^{2}-39189168 n+14092603\right) c_{n-1}&+ \\
20\left(4500 n^{2}-18900 n+19739\right) c_{n-2} + c_{n-3}  & = 0, 
\end{align*}
with initial terms $c_0 = 1, c_1 = -161/(2^{10}\cdot 3^{5})$ and $c_2 = 26605753/(2^{23} \cdot 3^{12} \cdot 5^2)$. Note that this recursion comes from a topological ODE in the sense of \cite{BeDuYa18}. Zagier mentions that $c_n$ decay like $1/n!^2$, therefore one may hope to find rational numbers $u, v$ and a non-zero integer $w$ such that\footnote{Recall that $(x)_n$  denotes the rising factorial: $(x)_n \coloneqq x\cdot (x+1)\cdots (x+n-1)$.}
\[
\widetilde{c}_n = c_n \cdot (u)_n \cdot (v)_n \cdot  w^n
\]
is an integer for every $n \geq 0$. Indeed, Zagier claims that using a formula by him and Yang, it can be shown that 
\[
a_n \coloneqq c_n \cdot (3/5)_n \cdot (4/5)_n \cdot  (2^{10} \cdot 3^5 \cdot 5^4)^n \in \ZZ,
\]
however he also mentions that they did not succeed in finding a closed expression for the generating function $\sum_{n\geq0} a_n x^n$, nor proving that it is an algebraic function~\cite[p.~769]{Zagier18}. 

As we will see, this is a perfect example for the power of \texttt{gfun}, using which we can easily answer questions about the sequence $(a_n)_{n\geq0}$ and its generating function. Moreover, this example demonstrates well the ``guess and prove'' strategy for P-recursive sequences.   

We start by defining 
\begin{align*}
&\texttt{> p0 := 80352000*n*(5*n - 1)*(5*n - 2)*(5*n - 4):}\\
&\texttt{> p1 := 25*(2592000*n\^{}4 - 16588800*n\^{}3 + 39118320*n\^{}2 }\\ 
& \qquad \qquad \qquad \texttt{- 39189168*n + 14092603):}\\
&\texttt{> p2 := 20*(4500*n\^{}2 - 18900*n + 19739):}
\end{align*}
such that the sequence $(c_n)_{n\geq0}$ is simply given by
\begin{align*}
    \texttt{> rec\_{}c := \{p0*c(n) + p1*c(n - 1) + p2*c(n - 2) + c(n - 3),}\\
    \texttt{c(0) = 1, c(1) = -161/(2\^{}10*3\^{}5), c(2) = 26605753/(2\^{}23*3\^{}12*5\^{}2)\}}:
\end{align*}
Using the \gfun function \texttt{rectoproc} we can convert the recurrent definition into a Maple procedure:
\[
\texttt{pro\_c := rectoproc(rec\_c, c(n)):}
\]
Hence the first few terms of the sequence $(a_n)_{n\geq0}$ can be computed easily: 
\begin{align*}
&\texttt{> w := 2\^{}10*3\^{}5*5\^{}4:}\\
&\texttt{> seq(pro\_c(n)*pochhammer(3/5, n)*pochhammer(4/5, n)*w\^{}n, n=0..3);}\\
&\quad 1, -48300, 7981725900, -1469166887370000
\end{align*}
\begin{enumerate}
    \item Use effective closure properties of P-recursive sequences.
    \item First guess and then prove the recursion. 
    \item Use the recurrence from $1.$ and find an equivalent one of minimal-order. 
    \item Guess and prove the differential equation for the generating function. Then convert it into a recurrence.
\end{enumerate}
We will now show that all four methods can be easily performed using Maple's \texttt{gfun}. We will see that they yield different (but correct) recursions/differential equations -- a fact which might look surprising at first glance.

\subsubsection{Effective closure properties}\label{sec:cp}
First we will find the recursion for $(a_n)_{n\geq0}$ using the closure properties implemented in \texttt{gfun}. We define the recurrences for the rising factorials:
\begin{align*}
&\texttt{> rec\_ph3 := \{(n + 3/5)*c(n) = c(n + 1), c(0) = 1\}:}\\
&\texttt{> rec\_ph4 := \{(n + 4/5)*c(n) = c(n + 1), c(0) = 1\}:}
\end{align*}
Now we compute the recurrences for $c_n \cdot (3/5)_n$, then for $c_n \cdot (3/5)_n \cdot (4/5)_n$ and finally for $a_n = c_n \cdot (3/5)_n \cdot (4/5)_n \cdot w^n$:
\begin{align*}
    &\texttt{> `rec*rec`(rec\_c, rec\_ph3, c(n)):}\\
    &\texttt{> `rec*rec`(\%, rec\_ph4, c(n)):}\\
    &\texttt{> rec\_a := `rec*rec`(\%, {c(0) = 1, c(n + 1) = c(n)*w}, c(n)):}
\end{align*}
This gives a proof that the sequence $(a_n)_{n\geq0}$ satisfies the recursion
\begin{align}\label{eq:reca}
    &p_3(n)a_{n+3} + p_2(n)a_{n+2}+p_1(n)a_{n+1}+p_0(n)a_{n} = 0, \quad \text{where}\\
    &p_3(n) = 31(n + 3)(5n + 11), \nonumber\\
    &p_2(n) = 60(2592000n^4 + 14515200n^3 + 29787120n^2 + 27559152n + 10644379), \nonumber\\
    &p_1(n) = 2^{14}\cdot3^6\cdot5^2 (5n + 8)(5n + 9)(4500n^2 + 8100n + 3539), \nonumber\\
    &p_0(n) = 2^{22}\cdot3^{11}\cdot5^3 (5n + 8)(5n + 3)(5n + 9)(5n + 4).  \nonumber\\\nonumber
\end{align}

\subsubsection{Guessing the recursion} \label{sec:guessrec}
A different way to find a recurrence relation for $a_n$ is to guess it first and then prove the guess. We first compute 51 terms of the recursion:
\begin{align*}
&\texttt{> a := [seq(pro\_c(n)*pochhammer(3/5,n)*pochhammer(4/5,n)*w\^{}n,}\\
& \qquad \qquad \texttt{ n = 0..50)]:}
\end{align*}
Then we use the function \texttt{listtorec} in order to guess a linear relation with polynomial coefficients:
\[
\texttt{> rec\_a\_guess := listtorec(a,u(n))[1]:}
\]
We find a much smaller recurrence of order 2 (compared to the one proven above of order 3):
\begin{align} \label{eq:guesseda}
    &\widetilde{p}_2(n)u_{n+2}+\widetilde{p}_1(n)u_{n+1}+\widetilde{p}_0(n)u_{n} = 0, \quad \text{where}\\
    &\widetilde{p}_2(n) = (5n + 6)(n + 2)(60n + 43), \nonumber\\
    &\widetilde{p}_1(n) = 300 (216000n^3 + 759600n^2 + 836940n + 290603)  \nonumber\\
    &\widetilde{p}_0(n) = 2^{12} \cdot 3^6 \cdot 5^2 (5n + 4)(5n + 3)(60n + 103). \nonumber
\end{align}
This recursion is found in a fraction of a second, however is not yet proven. Because we guessed it using 51 terms, we can only be certain that it gives correct terms $a_n$ for $0 \leq n \leq 50$. One can easily check by computing and comparing terms that this recurrence also holds true for $n \leq 100$ or $n \leq 1000$. 

There are several possibilities for proving the guess. Arguably the shortest one is explained in \S\ref{sec:minimalrec}. However, for pedagogical reasons, we will first argue on the level of differential operators, since this is exactly the procedure one would follow if trying to prove equality of two D-finite functions. First define $(\widetilde{a}_n)_{n\geq0}$ as being the unique sequence satisfying equation (\ref{eq:guesseda}) with initial terms $\widetilde{a}_0 = a_0$ and $\widetilde{a}_1 = {a}_1$.  Note that in order to guarantee uniqueness, we use that $\widetilde{p}_{2}(n)$ is non-zero for $n \geq 0$. Now we rigorously compute the differential equations satisfied by the generating functions of $({a}_n)_{n\geq0}$ and $(\widetilde{a}_n)_{n\geq0}$ using the \gfun function \texttt{rectodiffeq}:
\begin{align*}
    &\texttt{> deq\_a := rectodiffeq(rec\_a,c(n),y(x)):}\\
    &\texttt{> deq\_a\_guess := rectodiffeq(rec\_a\_guess,u(n),y(x)):}
\end{align*}
We find different differential equations of order 4 and 3 respectively. Now we translate both equations to differential operators using the Maple package \texttt{DEtools}.
\begin{align*}
    &\texttt{> L\_a := de2diffop(deq\_a[1],y(x),[Dx,x]):}\\
    &\texttt{> L\_a\_guess := de2diffop(deq\_a\_guess,y(x),[Dx,x]):}
\end{align*}
We can compute the LCLM $L$ of the two operators, rewrite it as a differential equation and transform it back to a recurrence for the coefficients of the solutions:
\begin{align*}
    &\texttt{> L := LCLM(L\_a, L\_a\_guess, [Dx,x]):}\\
    &\texttt{> deq := diffop2de(L, y(x), [Dx,x]):}\\
    &\texttt{> diffeqtorec(deq, y(x), u(n)):}
\end{align*}
We find exactly the same recurrence relation as in equation (\ref{eq:reca}). Notice that the leading coefficient $p_3(n)$ does not vanish for positive $n$. Therefore,  if the initial terms are prescribed to be $(a_0,a_1,a_2)$, the differential equation corresponding to the operator $L$ has the unique solution $\sum_{n\geq0} a_nx^n$. But since $L$ is defined as the LCLM of the operators corresponding to the sequences $(a_n)_{n\geq0}$ and $(\widetilde{a}_n)_{n\geq0}$, it also annihilates $\sum_{n\geq0} \widetilde{a}_nx^n$. This proves that $a_n  = \widetilde{a}_n$ for all $n \in \NN$ and consequently that our guessed recursion $(\ref{eq:guesseda})$ is correct.\\

\subsubsection{Minimal-order recursion} \label{sec:minimalrec}
With the newest Maple version of 2021 we have a great shortcut thanks to van Hoeij's improvement in the package \texttt{LREtools}. We can namely directly algorithmically find the minimal-order linear recurrence after obtaining the one in \S\ref{sec:cp} by just calling 
\[
    \texttt{> LREtools[`MinimalRecurrence`](rec\_a,u(n)):}
\]
We find exactly $(\ref{eq:guesseda})$. This not only yields another proof of the correctness of the guessed (and then proven) recurrence, but also proves its minimality. Note that this method does not rely on guessing. \\

\subsubsection{Guessing the ODE} \label{sec:guessdeq}
Finally, we can also guess the differential equation for $\sum_{n\geq0}a_nx^n$, prove its correctness and transform it to a recursion. This method has the advantage that we might discover a differential equation of smaller order than we would obtain by converting the recurrences above. We simply call the \gfun function
\[
    \texttt{> deq\_a\_ODEguess := listtodiffeq(a, y(x))[1]:}
\]
where \texttt{a} is the list of the first 51 terms of the sequence $(a_n)_{n\geq0}$ we computed in \S\ref{sec:guessrec}. We find a small differential equation of order 2 (compared to the differential equations above \texttt{deq\_a} and \texttt{deq\_a\_guess} of orders 4 and 3).
\begin{align} \label{eq:deqa}
    &q_2(x) y''(x) + q_1(x) y'(x) + q_0(x) y(x) = 0, \quad \text{where}\\
    &q_2(x) = 5x(302400x - 31)(373248000x^2 + 216000x + 1), \nonumber\\
    &q_1(x) = 1354442342400000x^3 + 64571904000x^2 - 61473600x - 31,\nonumber\\
    &q_0(x) = 300(902961561600x^2 - 240974784x - 4991).\nonumber
\end{align}
The proof of the correctness of this guess is similar to the proof in \S\ref{sec:guessrec}. In this case we actually found a (right) factor of \texttt{L\_a} as we can see by computing the GCRD (\texttt{L\_a\_guess2} is the corresponding differential operator to (\ref{eq:deqa})):
\[
    \texttt{> GCRD(L\_a, L\_a\_guess2, [Dx,x]):}
\]
This gives exactly \texttt{L\_a\_guess2}. Since our solutions to (\ref{eq:reca}) and (\ref{eq:deqa}) agree up to precision $3$, they must be equal. Therefore the guess must be correct.

Transforming this differential equation into a recursion for $(a_n)_{n\geq0}$ yields yet another recurrence, this time of order 3. This means that we found two different recurrences describing $(a_n)_{n\geq0}$ and three different differential equations describing the generating function. The orders of these objects are displayed in Table~\ref{tab:Zaga}.
\begin{table}[]
\centering
\begin{tabular}{l|c|c}
& \multicolumn{1}{l|}{Order of recurrence} & \multicolumn{1}{l}{Order of ODE} \\ \hline
Closure properties (\S\ref{sec:cp})     & 3  & 4 \\ \hline
\begin{tabular}[c]{@{}l@{}}Guessing the recurrence (\S\ref{sec:guessrec}) \\ Computing the minimal one (\S\ref{sec:minimalrec}) \end{tabular} & 2 & 3 \\ \hline
Guessing the ODE (\S\ref{sec:guessdeq})& 3 & 2
\end{tabular}
\caption{Orders of different recurrences and ODEs for the sequence $(a_n)_{n\geq0}$ and its generating function.}
\label{tab:Zaga}
\end{table}

\subsubsection{The generating function of $(a_n)_{n\geq0}$ is algebraic}
In this part we will prove that 
\begin{thm}\label{thm:aalg}
The generating function of the sequence $(a_n)_{n\geq0}$ is algebraic. 
\end{thm}
Using Maple we can actually solve the differential equation of order two for $\sum_{n\geq0} a_nx^n$ we found and proved in \S\ref{sec:guessdeq}. Simply calling
\[
\texttt{> dsolve(deq\_a\_ODEguess[1]):}
\]
shows that every solution of (\ref{eq:deqa}) is a linear combination of 
\begin{align*}
&A_1(x) \coloneqq u_1(x) \cdot \pFq{2}{1}{-1/60, 11/60}{2/3}{\frac{p_1(x)}{p_2(x)}} \quad \text{ and }\\
&A_2(x) \coloneqq u_2(x)\cdot \pFq{2}{1}{19/60, 31/60}{4/3}{\frac{p_1(x)}{p_2(x)}},
\end{align*}
where $u_1(x),u_2(x)$ are explicit algebraic functions and $p_1(x),p_2(x)$ are known polynomials and $\pFq{2}{1}{a, b}{c}{x}$ is the Gaussian hypergeometric function defined in~(\ref{eq:2F1}). 

Now we can proceed in two different ways. First, a classical work~\cite{Schwarz1872} by Schwarz from 1873, classifies all Gaussian hypergeometric functions that are algebraic. Applying this classification, known as \textit{Schwarz's list}, we can convince ourselves that both ${}_2F_1$'s above are algebraic. 
\begin{lem}\label{lem:f1f2}
The functions 
\[
f_1(x) \coloneqq \pFq{2}{1}{-1/60, 11/60}{2/3}{x} \text{ and } f_2(x) \coloneqq \pFq{2}{1}{19/60, 31/60}{4/3}{x}
\]
are algebraic.
\end{lem}

Another solution which is completely different in spirit, but useful also for more general problems of deciding algebraicity, is to use the ``guess and prove'' method explained in \S\ref{sec:guessnprove}. We can first try to guess and then to prove minimal polynomials for $f_1(x)$ and $f_2(x)$. This will not only provide a proof for algebraicity of the functions, but also give explicit minimal polynomials. In order to make the computations easier we will actually work with twelfth powers of $f_1(x)$ and $f_2(x)$. Moreover, here we only explain the computations for $f_1(x)$, because the exact same code works for $f_2(x)$ as well.

First we compute 100 terms of the series expansion for $f_1(x)^{12}$:
\begin{align*}
    &\texttt{> f12 := hypergeom([-1/60, 11/60], [2/3], x)\^{}12:}\\
    &\texttt{> ser1 := series(f12,x,100):} 
\end{align*}
Then we guess an annihilating polynomial for this series using \texttt{gfun}:
\begin{align*}
    &\texttt{> P := seriestoalgeq(ser1,y(x)):} 
\end{align*}
After a few seconds, this finds a polynomial $P(x,y) \in \ZZ[x,y]$ of degree $20$ in $y$ and $4$ in $x$. We note that $\partial_yP(0,1) \neq 0$, therefore there exists only one power series solution $f(x)$ to $P(x,y) = 0$ such that $f(0)=1$. Now we will confirm our guess. First we use the effective property that any algebraic function is D-finite:
\begin{align*}
    \texttt{> deq := algeqtodiffeq(P, y(x)):}
\end{align*}
Here we find an inhomogeneous differential equation, which we convert into a homogeneous one using \texttt{gfun}'s \texttt{diffeqtohomdiffeq}. Let us call the resulting equation \texttt{deqh}. It holds that any solution in $y(x)$ to $P(x,y)=0$ satisfies the differential equation \texttt{deqh}. Moreover, we can find the differential equation satisfied by $f_1(x)^{12}$ by simply calling
\[
    \texttt{> deqf12 := holexprtodiffeq(f12,y(x)):}
\]
We find exactly the same differential equation as \texttt{deqh}. By uniqueness and after checking enough terms, we can conclude that $P(x,y)$ indeed annihilates $f_{1}(x)^{12}$. Hence, $f_1(x)$ is algebraic. Moreover, the irreducible polynomial $P(x,y^{12})$ is then clearly the minimal polynomial for $f_{1}(x)$. This concludes the proof of Lemma~\ref{lem:f1f2} for the first function, while the second one can be done completely analogously. \\

Coming back to the generating function $\sum_{n\geq0} a_n x^n$, Lemma~\ref{lem:f1f2} implies that both $A_1(x)$ and $A_2(x)$ are algebraic. Then any linear combination of them must be an algebraic function as well. This proves Theorem~\ref{thm:aalg}.

\subsection{Monotonicity of \Iso}\label{sec:iso}
Our second example originates in biology, more specifically in the so-called Canham model which predicts the shape of biomembranes such as blood cells~\cite{Canham70}. Roughly speaking, the model asks to minimize the \emph{Willmore energy}
\[
W(S) = \int_S H^2 \mathrm{d}A,
\]
over orientable closed surfaces $S$ with prescribed genus, area and volume. The existence of a solution to Canham's model was investigated most notably by Schygulla~\cite{Schygulla12} and Keller, Mondino, Rivière~\cite{KeMoRi14}. We are rather interested in the uniqueness of the solution, studied by Seifert~\cite{Seifert97}, Chen et al.~\cite{ChYuBrKuYiZi21} and most recently in the article by Yu and Chen \cite{YuCh20}. 

Yu and Chen observed that the solution to this model is unique in the genus-one case if a certain function, called $\Iso(z)$, is strictly increasing for $z \in [0,\sqrt{2}-1)$. One way \cite[\S4.1]{YuCh20} to define $\Iso(z)$ is 
\begin{equation} \label{def:iso}
\Iso(z) \coloneqq {3 \cdot 2^{3/4} \pi^{-1/2}} \cdot \frac{\bar{V}(z^2)}{\bar{A}^{3/2}(z^2)},
\end{equation}
where $\bar{A}(z) = \sum_{n\geq0} a_n z^n \in \QQ\ps{z}$ and $\bar{V}(z) = \sum_{n\geq0} v_n z^n \in \QQ\ps{z}$ are given by 
\begin{align*}
    \bar{A}(z) &= \frac{1}{\sqrt{2}\pi^2}\int_0^{2\pi} \int_0^{2\pi} \frac{\sqrt{2} + \sin(v)}{Q(u,v,1\,;\sqrt{z})^2} \mathrm{d}u\mathrm{d}v, \text{ and}\\
    \bar{V}(z) &= \frac{1}{\sqrt{2}\pi^2}\int_0^{1} \int_0^{2\pi}\int_0^{2\pi} \frac{r\sqrt{2} + r^2\sin(v)}{Q(u,v,r\,;\sqrt{z})^3} \mathrm{d}u\mathrm{d}v\mathrm{d}r,
\end{align*}
where 
\[
Q(u,v,r\,;\,z) = {1+2(\sqrt{2}+r\sin(v))\cos(u)z+(2+r^2+2\sqrt{2}r\sin(v))z^2}.
\]
Less than one year later, $\Iso(z)$ was proven \cite{MeMe20} to be increasing:
\begin{thm} \label{thm:iso} 
    The function $\Iso(z)$ is strictly increasing on $[0,\sqrt{2}-1)$.
\end{thm}

Let us briefly explain the ideas that led to the proof of this theorem. Using the paradigm of creative telescoping, Yu and Chen found and proved linear recurrences of order $3$ with polynomial coefficients of degree $4$ for the sequences $(a_n)_{n\geq0}$ and $(v_n)_{n\geq0}$. We refer to Proposition 4.1 in \cite{YuCh20} for the explicit formulas. The same authors observed that $\Iso(z)$ is increasing if $\bar{D}(z) \geq 0$ on $z \in [0,\sqrt{2}-1)$, where 
\[
\bar{D}(z) \coloneqq 2 \bar{V}'(z)\bar{A}(z) - \bar{V}(z)\bar{A}'(z).
\]
Clearly, the function $\bar{D}(z) = \sum_{n\geq0} d_n z^n $ is D-finite, hence the sequence $(d_n)_{n\geq0}$ is  P-recursive. The positivity of the function obviously follows if $d_n \geq 0$ for all $n$. This is the central conjecture of the work by Yu and Chen and the main theorem in \cite{MeMe20} by Melczer and Mezzarobba. The latter authors used rigorous and effective analysis of the asymptotics of $(d_n)_{n\geq0}$ to show the positivity of this sequence and consequently settle Theorem~\ref{thm:iso}. 

Another proof of Theorem~\ref{thm:iso}, which is completely different in spirit, was recently proposed by Bostan and the author in \cite{BoYu21}. Even though one strength of the final version of this proof is that it is completely elementary and may be verified without a computer, one should admit that finding this solution required algorithms from \texttt{gfun} and Maple in general. Exploring the discovery of this proof of monotonicity of $\Iso(z)$ will be our second example. \\

The first natural idea to show that $\Iso(z)$ is strictly increasing is to find an explicit closed formula for it. Note that the function is not D-finite, but given as the ratio of two holonomic functions. Therefore we will try to find explicit formulas for the numerator and denominator separately.

Converting the recurrence relations for $(a_n)_{n\geq0}$ and $(v_n)_{n\geq0}$ into differential equations for $\bar{A}(z)$ and $\bar{V}(z)$ using \texttt{gfun}'s \texttt{rectodiffeq}, we find two linear differential equations of order 4 and degree 5. First, it seems that working with these differential operators is hopeless, however it turns out that they are not minimal for our functions of interest. In fact, we can easily guess and then prove the much smaller minimal differential equations for $\bar{A}(z)$ and $\bar{V}(z)$: we first compute two lists of 100 terms of $(a_n)_{n\geq0}$ and $(v_n)_{n\geq0}$, call them \texttt{la} and \texttt{lv}. Then in less than one second we can execute
\begin{align*}
    &\texttt{> deqA := listtodiffeq(la,y(x));}\\
    &\texttt{> deqV := listtodiffeq(lv,y(x));}
\end{align*}
This finds two second-order differential equations for $\bar{A}(z)$ and $\bar{V}(z)$. Similarly to the example in \S\ref{sec:guessdeq}, it is not difficult to prove the correctness of both guesses by computing the GCRD of the known operators with the new ones and arguing by uniqueness of solutions. 

Now it turns out that we can actually solve the new smaller differential equations using Maple's \texttt{dsolve}. We find that 
\begin{equation}\label{eq:A}
\bar{A}(z) =   
\frac{4(z + 1)}{(1 - 6z + z^2)^{3/2}}  
\, \cdot \, \pFq{2}{1}{-\frac{1}{2},\frac{3}{2}}{1}{\frac{-4 z}{1-6z+z^2}},
\end{equation}
and
\begin{equation}\label{eq:V}
\bar{V}(z) =
\frac{2}{(1 - 6z + z^2)^{3/2}}  
\, \cdot \, \pFq{2}{1}{-\frac{3}{2},\frac{5}{2}}{1}{\frac{-4 z}{1 - 6z + z^2}}.
\end{equation}
As soon as these formulas for $\bar{A}(z)$ and $\bar{V}(z)$ are discovered, they can be verified ``by hand'' without sophisticated algorithms. In other words, \cite[Thm.~1]{BoYu21} has a human proof, but admittedly it had a computer-assisted discovery. 

The example, however, does not end here: we still need to prove that $\Iso(z) = {3 \cdot 2^{3/4} \pi^{-1/2}} \cdot {\bar{V}(z^2)} \cdot {\bar{A}^{-3/2}(z^2)}$ is monotonic, and here we will again use ``guess-and-prove'' and \texttt{gfun}. After a few simple manipulations \cite[p.~3]{BoYu21} of the definition of $\Iso(z)$, equations (\ref{eq:A}) and (\ref{eq:V}), and Gauss' summation theorem, one finds
\[
\Iso(z) = \frac{3}{2^{5/4}\cdot \sqrt{\pi}} \cdot w_{1/2}(x)^{-3/2} \cdot w_{3/2}(x),
\]
where $x=4z^2/(1-z^2)^2$ and $w_a(x) \coloneqq \pFq{2}{1}{-a,-a}{1}{x}\cdot (1+x)^{-a}$. The monotonicity of $\Iso$ follows if we can show that $w_a$ is decreasing on $[0,1]$ for $0<a<1$ and increasing on this interval if $a>1$. After calculating the derivative of $w_a(x)$, clearing the denominator and normalizing for unit constant coefficient, we find
\[
\frac{w_a'(x) }{a(a-1)} \cdot (1+x)^{a+1} = 1 + \frac{1}{2}(a - 3)ax + \frac{1}{12}(a - 5)(a - 1)(a - 2)ax^2 + O(x^3).
\]
The denominator above is $a(a-1)$ and therefore conveniently takes care of the case distinction in $0<a<1$ and $a>1$. Now we just need to argue that $w_a'(x) \cdot {(x+1)^{a+1}}/(a(a-1))$ is positive on $[0,1]$. This would be obvious if all Taylor coefficients of the function were positive. Unfortunately, this does not hold. However, after a few tries we find that
\[
    h(x) \coloneqq \frac{w_a'(x) }{a(a-1)} \cdot \frac{(1+x)^{a+1}}{(1-x)^{2a}} = 1 + \frac{1}{2}(a + 1)ax + \frac{1}{12}(a + 2)(a + 1)^2ax^2 + O(x^3)
\]
seems to have positive Taylor coefficients. If we calculate the series expansion of $h(x)$ to degree 10 and use \texttt{seriestorec}, we can guess an easy first-order recurrence relation for the coefficients of $h(x)$:
\[
(n + 2)(n + 1)u_{n+1} - (a + n + 1)(a + n)u_n = 0.
\]
Together with the initial term $u_0=1$ this (guessed) recurrence implies that 
\begin{equation}\label{eq:hgauss}
    \frac{w_a'(x) }{a(a-1)} \cdot \frac{(1+x)^{a+1}}{(1-x)^{2a}} = \pFq{2}{1}{a,a+1}{2}{x}.
\end{equation}
The right-hand side obviously has positive Taylor coefficients, hence the left-hand side is also positive on $[0,1)$. Therefore, modulo the proof of equation~(\ref{eq:hgauss}), this proves Theorem~\ref{thm:iso}.

A ``human'' proof of identity~(\ref{eq:hgauss}) is the content of Lemma~1 in \cite{BoYu21} which relies on Gauss' contiguous relations. Here we will present another algorithmic proof, which exploits the fact that both sides of the identity are D-finite. 

For the left-hand side of (\ref{eq:hgauss}) (stored as \texttt{h} in Maple) we find a differential equation and the corresponding operator by simply calling 
\begin{align*}
    &\texttt{> deq := holexprtodiffeq(h,y(x)):}   \\ &\texttt{> L := de2diffop(deq[1],y(x),[Dx,x]):}
\end{align*}
This gives a differential equation of order 4. The right-hand side of (\ref{eq:hgauss}) satisfies the second-order differential equation
\[
a(a + 1)y(x) + 2(ax + x - 1)y'(x) + x(x - 1)y''(x) = 0,
\]
and we will call the corresponding operator \texttt{L\_guess}. The output of 
\[
    \texttt{> GCRD(L,L\_guess,[Dx,x]);}
\]
is exactly \texttt{L\_guess}. Therefore, \texttt{L\_guess} (right-)divides \texttt{L}. Moreover, since the leading term of the recursion corresponding to \texttt{deq} is $(n + 5)(n + 3)(n + 4)^2(a + 3)$, the equation $\texttt{L}y=0$ has a unique solution if the first four Taylor coefficients of $y$ are prescribed. This is easily checked and hence (\ref{eq:hgauss}) is proved.\\

\bigskip

The two examples demonstrate several things at the same time. First, we saw that guessing a minimal-order operator is quite simple with \gfun and that the proof is always an easy argument on the level of differential operators. Moreover, both examples show in different ways that having explicit solutions in terms of Gaussian hypergeometric functions can be very useful in practice. The example \S\ref{sec:zag} also demonstrates the ``guess and prove'' strategy for proving algebraicity of a generating function. In the second example \S\ref{sec:iso} we saw that guessing can be used not only to predict identities like (\ref{eq:hgauss}) but also to simplify expressions like $h(x)$. Finally, the second application also demonstrates that both \gfun and \texttt{DEtools} work if parameters are involved as well.\\

\noindent \textbf{Acknowledgments:} The author is deeply indebted to Alin Bostan who was not only a great source of motivation for this article, but is also a dedicated teacher and enthusiast of most techniques and ideas explained in it. Moreover, the author is also grateful to Bruno Salvy for the clarifications of several questions regarding the new version of \texttt{gfun}. We also wish to thank Alin Bostan, Giancarlo Castellano and the anonymous referees for careful reading of the manuscript and constructive comments. Finally, we thank Robert Corless for pointing out some mistakes in the almost final version of the manuscript.

This work was financially supported by the DOC fellowship of the \href{https://www.oeaw.ac.at/en/}{ÖAW} (26101), the WTZ collaboration project of the \href{https://oead.at/en/}{OeAD} (FR 09/2021) and the {\href{https://specfun.inria.fr/chyzak/DeRerumNatura/}{DeRerumNatura}} project ANR-19-CE40-0018.

\bibliographystyle{alphaabbr}
\bibliography{bib_maple}

\newcommand{\etalchar}[1]{$^{#1}$}
\begin{thebibliography}{BCDVW20}

\bibitem[Ap{\' e}79]{Apery79}
R.~Ap{\' e}ry.
\newblock Irrationalit\'{e} de {$\zeta(2)$} et {$\zeta(3)$}.
\newblock {\em Ast\'{e}risque}, (61):11--13, 1979.

\bibitem[AS64]{AbSt64}
M.~Abramowitz and I.~A. Stegun.
\newblock {\em Handbook of Mathematical Functions with Formulas, Graphs, and
  Mathematical Tables}.
\newblock Dover, New York, ninth dover printing, tenth gpo printing edition,
  1964.

\bibitem[BCDVW20]{BaClDiWe20}
M.~Barkatou, T.~Cluzeau, L.~Di~Vizio, and J.-A. Weil.
\newblock Reduced forms of linear differential systems and the intrinsic
  {G}alois-{L}ie algebra of {K}atz.
\newblock {\em SIGMA Symmetry Integrability Geom. Methods Appl.}, 16:Paper No.
  054, 13, 2020.

\bibitem[BDY18]{BeDuYa18}
M.~Bertola, B.~Dubrovin, and D.~Yang.
\newblock Simple {L}ie algebras and topological {ODE}s.
\newblock {\em Int. Math. Res. Not. IMRN}, (5):1368--1410, 2018.

\bibitem[BK10]{BoKa10}
A.~Bostan and M.~Kauers.
\newblock The complete generating function for {G}essel walks is algebraic.
\newblock {\em Proc. Amer. Math. Soc.}, 138(9):3063--3078, 2010.
\newblock With an appendix by Mark van Hoeij.

\bibitem[BL97]{BeLa97}
B.~Beckermann and G.~Labahn.
\newblock Recursiveness in matrix rational interpolation problems.
\newblock volume~77, pages 5--34. 1997.
\newblock ROLLS Symposium (Leipzig, 1996).

\bibitem[BLS17]{BoLaSa17}
A.~Bostan, P.~Lairez, and B.~Salvy.
\newblock Multiple binomial sums.
\newblock {\em J. Symbolic Comput.}, 80(part 2):351--386, 2017.

\bibitem[BMP00]{BoPe00}
M.~Bousquet-M\'{e}lou and M.~Petkov\v{s}ek.
\newblock Linear recurrences with constant coefficients: the multivariate case.
\newblock volume 225, pages 51--75. 2000.
\newblock Formal power series and algebraic combinatorics (Toronto, ON, 1998).

\bibitem[BMP03]{BoPe03}
M.~Bousquet-M\'{e}lou and M.~Petkov\v{s}ek.
\newblock Walks confined in a quadrant are not always D-Finite.
\newblock {\em Theor. Comput. Sci.}, 307(2):257–276, October 2003.

\bibitem[Bos17]{Bostan17}
A.~Bostan.
\newblock {\em Computer algebra for lattice path combinatorics}.
\newblock {HDR} (accreditation to supervise research), Univ. Paris 13, 2017.
\newblock 79 pages.

\bibitem[Bos21]{Bostan21}
A.~Bostan.
\newblock Computer algebra in the service of enumerative combinatorics.
\newblock In {\em I{SSAC} '21---{P}roceedings of the 2021 {I}nternational
  {S}ymposium on {S}ymbolic and {A}lgebraic {C}omputation}, pages 1--8. ACM,
  New York, [2021] \copyright 2021.

\bibitem[BRS21]{BoRiSa21}
A.~Bostan, T.~Rivoal, and B.~Salvy.
\newblock Explicit degree bounds for right factors of linear differential
  operators.
\newblock {\em Bull. Lond. Math. Soc.}, 53(1):53--62, 2021.

\bibitem[BY22]{BoYu21}
A.~Bostan and S.~Yurkevich.
\newblock A hypergeometric proof that {$\sf{Iso}$} is bijective.
\newblock {\em Proc. Amer. Math. Soc.}, 150(5):2131--2136, 2022.

\bibitem[Can70]{Canham70}
P.~Canham.
\newblock The minimum energy of bending as a possible explanation of the
  biconcave shape of the human red blood cell.
\newblock {\em Journal of Theoretical Biology}, 26(1):61--81, 1970.

\bibitem[CC86]{ChCh86}
D.~V. Chudnovsky and G.~V. Chudnovsky.
\newblock On expansion of algebraic functions in power and {P}uiseux series.
  {I}.
\newblock {\em J. Complexity}, 2(4):271--294, 1986.

\bibitem[CC90]{ChCh90}
D.~V. Chudnovsky and G.~V. Chudnovsky.
\newblock Computer algebra in the service of mathematical physics and number
  theory.
\newblock In {\em Computers in mathematics ({S}tanford, {CA}, 1986)}, volume
  125 of {\em Lecture Notes in Pure and Appl. Math.}, pages 109--232. Dekker,
  New York, 1990.

\bibitem[Chy14]{Chyzak14}
F.~Chyzak.
\newblock {\em The {ABC} of Creative Telescoping: Algorithms, Bounds,
  Complexity}.
\newblock {HDR} (accreditation to supervise research), University Paris-Sud 11,
  April 2014.
\newblock 64 pages.

\bibitem[CK17]{ChKa17}
S.~Chen and M.~Kauers.
\newblock Some open problems related to creative telescoping.
\newblock {\em J. Syst. Sci. Complex.}, 30(1):154--172, 2017.

\bibitem[Com64]{Comtet64}
L.~Comtet.
\newblock Calcul pratique des coefficients de {T}aylor d'une fonction
  alg\'{e}brique.
\newblock {\em Enseign. Math. (2)}, 10:267--270, 1964.

\bibitem[CYB{\etalchar{+}}21]{ChYuBrKuYiZi21}
J.~Chen, T.~Yu, P.~Brogan, R.~Kusner, Y.~Yang, and A.~Zigerelli.
\newblock Numerical methods for biomembranes: conforming subdivision methods
  versus non-conforming {PL} methods.
\newblock {\em Math. Comp.}, 90(328):471--516, 2021.

\bibitem[Gra00]{Gray00}
J.~J. Gray.
\newblock {\em Linear differential equations and group theory from {R}iemann to
  {P}oincar\'{e}}.
\newblock Birkh\"{a}user Boston, Inc., Boston, MA, second edition, 2000.

\bibitem[IvH15]{ImHo15}
E.~Imamoglu and M.~van Hoeij.
\newblock Computing Hypergeometric Solutions of Second Order Linear
  Differential Equations Using Quotients of Formal Solutions.
\newblock In {\em Proceedings of the 2015 ACM on International Symposium on
  Symbolic and Algebraic Computation}, ISSAC '15, page 235–242, New York, NY,
  USA, 2015. Association for Computing Machinery.

\bibitem[Jun31]{Jungen31}
R.~Jungen.
\newblock Sur les s\'{e}ries de {T}aylor n'ayant que des singularit\'{e}s
  alg\'{e}brico-logarithmiques sur leur cercle de convergence.
\newblock {\em Comment. Math. Helv.}, 3(1):266--306, 1931.

\bibitem[KMR14]{KeMoRi14}
L.~G.~A. Keller, A.~Mondino, and T.~Rivi\`ere.
\newblock Embedded surfaces of arbitrary genus minimizing the {W}illmore energy
  under isoperimetric constraint.
\newblock {\em Arch. Ration. Mech. Anal.}, 212(2):645--682, 2014.

\bibitem[Kra99]{Krattenthaler99}
C.~Krattenthaler.
\newblock Advanced determinant calculus.
\newblock volume~42, pages Art. B42q, 67. 1999.
\newblock The Andrews Festschrift (Maratea, 1998).

\bibitem[KvH13]{KuHo13}
V.~J. Kunwar and M.~van Hoeij.
\newblock Second order differential equations with hypergeometric solutions of
  degree three.
\newblock In {\em I{SSAC} 2013---{P}roceedings of the 38th {I}nternational
  {S}ymposium on {S}ymbolic and {A}lgebraic {C}omputation}, pages 235--242.
  ACM, New York, 2013.

\bibitem[Mez10]{Mezzarobba10}
M.~Mezzarobba.
\newblock Num{G}fun: a package for numerical and analytic computation and
  {D}-finite functions.
\newblock In {\em I{SSAC} 2010---{P}roceedings of the 2010 {I}nternational
  {S}ymposium on {S}ymbolic and {A}lgebraic {C}omputation}, pages 139--146.
  ACM, New York, 2010.

\bibitem[MM20]{MeMe20}
S.~Melczer and M.~Mezzarobba.
\newblock {Sequence Positivity Through Numeric Analytic Continuation:
  Uniqueness of the {C}anham Model for Biomembranes}, 2020.
\newblock Technical Report
  \href{https://arxiv.org/abs/2011.08155}{arXiv:2011.08155} [math.CO].

\bibitem[Ore32]{Ore32}
O.~Ore.
\newblock Formale {T}heorie der linearen {D}ifferentialgleichungen. ({E}rster
  {T}eil).
\newblock {\em J. Reine Angew. Math.}, 167:221--234, 1932.

\bibitem[Poo60]{Poole60}
E.~G.~C. Poole.
\newblock {\em Introduction to the theory of linear differential equations}.
\newblock Dover Publications, Inc., New York, 1960.

\bibitem[Pó78]{Polya78}
G.~Pólya.
\newblock Guessing and Proving.
\newblock {\em The Two-Year College Mathematics Journal}, 9(1):21--27, 1978.

\bibitem[Sal05]{Salvy05}
B.~Salvy.
\newblock D-finiteness: algorithms and applications, 2005.
\newblock Invited talk in {\em Proceedings of ISSAC'05}.

\bibitem[Sal19]{Salvy19}
B.~Salvy.
\newblock Linear differential equations as a data structure.
\newblock {\em Found. Comput. Math.}, 19(5):1071--1112, 2019.

\bibitem[Sch73]{Schwarz1872}
H.~A. Schwarz.
\newblock \"{U}ber diejenigen {F}\"{a}lle, in welchen die {G}au{\ss}ische
  hypergeometrische {R}eihe einer algebraische {F}unktion ihres vierten
  {E}lementes darstellt.
\newblock {\em J. Reine Angew. Math.}, 75:292--335, 1873.

\bibitem[Sch12]{Schygulla12}
J.~Schygulla.
\newblock Willmore minimizers with prescribed isoperimetric ratio.
\newblock {\em Arch. Ration. Mech. Anal.}, 203(3):901--941, 2012.

\bibitem[Sei97]{Seifert97}
U.~Seifert.
\newblock Configurations of fluid membranes and vesicles.
\newblock {\em Advances in Physics}, 46(1):13--137, 1997.

\bibitem[ST20]{OEIS}
N.~J.~A. Sloane and {The OEIS Foundation Inc}.
\newblock The on-line encyclopedia of integer sequences.
\newblock \href{http://oeis.org/}{http://oeis.org/}, 2020.

\bibitem[Sta80]{Stanley80}
R.~P. Stanley.
\newblock Differentiably finite power series.
\newblock {\em European J. Combin.}, 1(2):175--188, 1980.

\bibitem[SZ94]{SaZi94}
B.~Salvy and P.~Zimmermann.
\newblock GFUN: A Maple Package for the Manipulation of Generating and
  Holonomic Functions in One Variable.
\newblock {\em ACM Trans. Math. Softw.}, 20(2):163–177, June 1994.

\bibitem[vdH99]{Hoeven99}
J.~van~der Hoeven.
\newblock Fast evaluation of holonomic functions.
\newblock {\em Theoret. Comput. Sci.}, 210(1):199--215, 1999.

\bibitem[vdH01]{Hoeven01}
J.~van~der Hoeven.
\newblock Fast evaluation of holonomic functions near and in regular
  singularities.
\newblock {\em J. Symbolic Comput.}, 31(6):717--743, 2001.

\bibitem[vdH07]{Hoeven07}
J.~van~der Hoeven.
\newblock Efficient accelero-summation of holonomic functions.
\newblock {\em J. Symbolic Comput.}, 42(4):389--428, 2007.

\bibitem[vdP79]{vdPoorten79}
A.~van~der Poorten.
\newblock A proof that {E}uler missed{$\ldots $}{A}p\'{e}ry's proof of the
  irrationality of {$\zeta (3)$}.
\newblock {\em Math. Intelligencer}, 1(4):195--203, 1978/79.
\newblock An informal report.

\bibitem[vdPS03]{PuSi03}
M.~van~der Put and M.~F. Singer.
\newblock {\em Galois theory of linear differential equations}, volume 328 of
  {\em Grundlehren der mathematischen Wissenschaften [Fundamental Principles of
  Mathematical Sciences]}.
\newblock Springer-Verlag, Berlin, 2003.

\bibitem[vHV15]{HoVi15}
M.~van Hoeij and R.~Vid\={u}nas.
\newblock Belyi functions for hyperbolic hypergeometric-to-{H}eun
  transformations.
\newblock {\em J. Algebra}, 441:609--659, 2015.

\bibitem[YC22]{YuCh20}
T.~Yu and J.~Chen.
\newblock Uniqueness of {C}lifford torus with prescribed isoperimetric ratio.
\newblock {\em Proc. Amer. Math. Soc.}, 150(4):1749--1765, 2022.

\bibitem[Zag18]{Zagier18}
D.~Zagier.
\newblock The arithmetic and topology of differential equations.
\newblock In {\em European {C}ongress of {M}athematics}, pages 717--776. Eur.
  Math. Soc., Z\"{u}rich, 2018.

\bibitem[Zei91]{Zeilberger91}
D.~Zeilberger.
\newblock The method of creative telescoping.
\newblock {\em J. Symbolic Comput.}, 11(3):195--204, 1991.

\end{thebibliography}

\end{document}